\documentclass[12pt]{amsart}
\usepackage{amsfonts}

\theoremstyle{definition}

\theoremstyle{remark}

\newcommand{\const}{\mathop{\rm const}\limits}

\newcommand{\supp}{\mathop{\rm supp}\limits}

\begin{document}

\begin{center}

{\bf LEBESGUE-RIESZ NORM ESTIMATES   \\

\vspace{3mm}

FOR FRACTIONAL LAPLACE TRANSFORM} \par

\vspace{3mm}

{\bf E. Ostrovsky}\\

e-mail: eugostrovskyo@list.ru \\

\vspace{3mm}

{\bf L. Sirota}\\

e-mail: sirota3@bezeqint.net  \\

\vspace{3mm}

Department of Mathematics and Statistics, Bar-Ilan University, 59200, Ramat Gan, Israel. \\

\vspace{3mm}

\begin{center}

   Abstract. \\

\end{center}

{\it We obtain in this short article the bilateral non-asymptotic estimations for the norm in Lebesgue-Riesz and
bilateral Grand Lebesgue spaces of the so-called fractional Laplace integral transform.\par

 We give also examples to show the sharpness of these inequalities.} \\

\end{center}

\vspace{3mm}


\vspace{3mm}


\vspace{3mm}

\section{Introduction. Notations. Statement of problem.}

\vspace{3mm}

 Let $ f: (0,\infty) \to R  $ be measurable function  and $ \kappa, r $  be constant real numbers.
The following linear operator (transform) $ L_{\kappa, r}[f](s) $ is named in the article \cite{Treuman1} as  a "Fractional Laplace Transform (FLT)":

$$
L_{\kappa, r}[f](s) \stackrel{def}{=} \int_0^{\infty} (1 + s t/\kappa)^{-\kappa - r} \ f(t) \ dt, \ s \ge 0. \eqno(1.1)
$$

 The classical Laplace transform $ L[f](s)  $ may be obtained formally as a limit

 $$
 L[f](s)  = \lim_{\kappa \to \infty} L_{\kappa, r}[f](s) = \int_0^{\infty} e^{-st} \ f(t) \ dt. \eqno(1.2)
 $$

 The applications of these transform in the statistical mechanics, differential equations, geophysics  etc. are described in
 the articles  \cite{Rahimy1}, \cite{Sokolov1}, \cite{Treuman0}, \cite{Treuman1}.\par

 Many properties of FLT: inversion, convolution identity, transform of derivatives etc. are investigated in the article
\cite{Treuman1}.\par

\vspace{3mm}

 {\bf  We intent to obtain in this report the $ L_p \to L_q  $ operator norm estimates for the fractional Laplace operator of the form }

 $$
 | \ L_{\kappa, r}[f](\cdot ) \ |_q \le K_{\kappa,r}(p,q) \ |f|_p, \ p,q \ge 1. \eqno(1.3)
 $$

\vspace{3mm}

 Hereafter

$$
|f|_p := \left[ \int_0^{\infty} |f(t) |^p \ dt \right]^{1/p};
$$
the {\it weight} case will be consider further.\par

We agree to take as a capacity to the number $ K_{\kappa,r}(p,q) $ its minimal value:

$$
K_{\kappa,r}(p,q) \stackrel{def}{=} \sup_{f: |f|_p \in (0,\infty)} \left[ \frac{| \ L_{\kappa, r}[f](\cdot ) \ |_q}{|f|_p}   \right]. \eqno(1.4)
$$

 The estimates of a form (1.3) for "pure" Laplace transform is described in \cite{Hardy1}, \cite{Okikiolu1},
\cite{Ostrovsky100}, \cite{Ostrovsky102}.\par

\vspace{4mm}

 We use symbols $C(X,Y),$ $C(p,q;\psi),$ etc., to denote positive
constants along with parameters they depend on, or at least
dependence on which is essential in our study. To distinguish
between two different constants depending on the same parameters
we will additionally enumerate them, like $C_1(X,Y)$ and
$C_2(X,Y).$ The relation $ g(\cdot) \asymp h(\cdot), \ p \in (A,B), $
where $ g = g(p), \ h = h(p), \ g,h: (A,B) \to R_+, $
denotes as usually

$$
0< \inf_{p\in (A,B)} h(p)/g(p) \le \sup_{p \in(A,B)}h(p)/g(p)<\infty.
$$
The symbol $ \sim $ will denote usual equivalence in the limit
sense.\par
We will denote as ordinary the indicator function
$$
I(x \in A) = 1, x \in A, \ I(x \in A) = 0, x \notin A;
$$
here $ A $ is a measurable set.\par
 All the passing to the limit in this article may be grounded by means
 of Lebesgue dominated convergence theorem.\par

\bigskip

\section{Main result: upper estimations for fractional Laplace operator}

\vspace{3mm}

 Some notations. Let $ p = \const \ge 1; $ we denote as ordinary by $ p' $ its conjugate number: $ p' = p/(p-1),  $ where
$ 1' := + \infty. $ \par
 We impose the following restriction on the values $ (p,q) $ in the inequality (1.3):

$$
q = p' = p/(p-1) \hspace{6mm} \Leftrightarrow \frac{1}{p} +  \frac{1}{q} = 1. \eqno(2.1)
$$

 Denote also  for the values $ (\kappa,r): \ \kappa + r > 1/2 $

$$
v(-1/2, \kappa + r) = \int_0^{\infty} \frac{x^{-1/2}}{(1+ x)^{\kappa + r}} \ dx =
\sqrt{\pi} \cdot \frac{\Gamma(\kappa + r - 1/2)}{\Gamma(\kappa + r)}, \eqno(2.2)
$$

$$
w(\kappa,r) = \sqrt{\kappa} \ v(-1/2, \kappa + r) = \sqrt{\pi \ \kappa} \cdot \frac{\Gamma(\kappa + r - 1/2)}{\Gamma(\kappa + r)},
$$

$$
z_{\kappa,r}(p) = w(\kappa,r)^{2/p'} = w(\kappa,r)^{2/q}. \eqno(2.3)
$$

\vspace{3mm}

{\bf Theorem 2.1.} {\it  Suppose in addition    }

$$
 p \in [1,2], \  \kappa + r > 1/2.  \eqno(2.4)
$$
{\it Then we deduce the following upper estimate: }

$$
K_{\kappa,r}(p,q) \le z_{\kappa,r}(p), \eqno(2.5a)
$$

{\it or, in detail:  }

$$
|L_{\kappa,r}[f]|_q =  |L_{\kappa,r}[f]|_{p'} \le
\left[ \sqrt{\pi \ \kappa} \cdot \frac{\Gamma(\kappa + r - 1/2)}{\Gamma(\kappa + r)} \right]^{2/p'} \cdot |f|_p. \eqno(2.5b)
$$

\vspace{3mm}

{\bf Proof} is at the same as one for the Laplace transform, see, e.g.  \cite{Okikiolu1}, p. 341-350. Namely,
let us consider the following linear operator:

$$
L_0 V_p [f](s) = \int_0^{\infty}  t^{2/p - 2} \ h(s/t) \ f(t) \ dt,
$$

$$
h(x) = h_{\kappa,r}(x) = ( 1 + x/\kappa)^{- \kappa - r}.
$$

 It is proved in fact in \cite{Okikiolu1}, p. 346

$$
|L_{\kappa,r}[f]|_q  = |L_0 V_p [f] |_q \le  \left[  \int_0^{\infty} t^{-1/2} \ h(t) \ dt  \right]^{2/p'}   \cdot  |f|_p =
z_{\kappa,r}(p) \cdot |f|_p,
$$
Q.E.D.\\

\vspace{3mm}

 {\bf Remark 2.1.  } Note that we derive as $ \kappa \to \infty $  from the inequality (2.5b) the $ L_p-L_q $ estimation for "pure"
 Laplace transform. \par

\bigskip

\section{Main result: inverse  estimations for fractional Laplace operator}

\vspace{3mm}

{\it A. \  Necessity of our conditions.}\\

\vspace{3mm}

 The necessity of restrictions (2.4) is proved in \cite{Okikiolu1}, p. 220-222; we need to prove the relation (2.1). \par

\vspace{3mm}

{\bf Theorem 3.1.} {\it   Suppose for some values $ (p,q), \ p,q > 0 $
there exists a finite constant (function) $ K = K_{\kappa,r}(p,q) $ such that
the inequality (1.3) is satisfied for arbitrary function $ f $  from the Schwartz class $ C_0^{\infty}(R_+). $
Then the relation (2.1) there holds. }\par

\vspace{3mm}

{\bf Proof.} We will apply the well known scaling method; see, e.g.   \cite{Stein1}, p. 285;  \cite{Talenti1}, \cite{Ostrovsky100}, \cite{Ostrovsky101}, \cite{Ostrovsky102}. Indeed, let $ f \ne 0 $ be any function from the set  $ C_0^{\infty}(R_+), $  then

$$
|y|_q \le K_{\kappa,r}(p,q) \ |f|_p, \ y:= L_{\kappa,r}[f]. \eqno(3.1)
$$

 Let $  \lambda  $ be arbitrary positive number. Introduce the so-called dilation operator, more exactly, the family of operators
$ T_{\lambda}: $

$$
T_{\lambda}[f](x) = f(\lambda \cdot x).
$$
 Obviously, $ T_{\lambda}[f] \in C_0^{\infty}(R_+). $ We derive substituting into (3.1) the function  $ T_{\lambda}[f]  $
instead  the initial function $ f: $

$$
|y_{\lambda}|_q \le K_{\kappa,r}(p,q) \ |T_{\lambda}[f]|_p, \ y:= L_{\kappa,r}[T_{\lambda}f]. \eqno(3.2)
$$
 But

$$
| T_{\lambda} \ f|_p = \lambda^{-1/p} \ |f|_p,  \hspace{6mm} |y_{\lambda}|_q = \lambda^{-1 + 1/q} \ |L_{\kappa,r} \ [f]|_q  \le
\lambda^{-1 + 1/q} \ K |f|_p,
$$

Therefore

$$
\lambda^{-1 + 1/q} \ K |f|_p \ge \lambda^{-1/p} |f|_p. \eqno(3.3)
$$
 Since the number $ \lambda  $ is arbitrary positive, we deduce from (3.3)

$$
-1 + 1/q = - 1/p,  \ \Leftrightarrow 1/p + 1/q = 1.
$$

\vspace{3mm}

 {\bf Remark 3.1.} The proposition of theorem 3.1 remains true, with at the same proof,
for the more general integral transforms of the form

$$
y(s) = \int_{R} h(s \cdot t) \ f(t) \ dt.
$$

\vspace{3mm}

{\it B. \ Lower estimates. } \\

  Note that the unique critical point of the values $  p  $ is the value $  p = 2-0. $  We can restrict ourselves the values
  $  p  $ by the set $ 1 \le p < 2. $  \par

 Let us choose as a capacity of trial function

$$
f_0(x) := x^{-1/2} \ I_{(0,1)}(x),
$$
then

$$
|f_0|_p^p = \int_0^1 x^{-p/2} \ dx = \frac{2}{2-p}; \ |f_0|_p = \left[ \frac{2}{2-p} \right]^{1/p} \asymp  \left[ \frac{2}{2-p} \right]^{1/2};
$$

$$
g_0(s):= L_{\kappa,r} [f](s) =\int_0^1 (1 + sx/\kappa)^{- \kappa - r  } \ \frac{dx}{\sqrt{x}}= s^{-1/2}
\int_0^s (1 + z/\kappa)^{-\kappa - r} \frac{dz}{\sqrt{z}}
$$
 and we have as $  s \to \infty: $

$$
g_0(s) \sim s^{-1/2} \ \int_0^{\infty} (1 + z/\kappa)^{- \kappa - r} \ \frac{dz}{\sqrt{z}}=
$$

$$
\sqrt{ \kappa/s} \  v(-1/2, \kappa + r) =  \sqrt{\pi \ \kappa/s } \cdot \frac{\Gamma(\kappa + r - 1/2)}{\Gamma(\kappa + r)}.
$$

 Therefore, we have as $ p \to 2-0  $ or equally as $ q \to 2 +0: $

$$
|L_{\kappa,r}[f_0]|_q =  |L_{\kappa,r}[f_0]|_{p'} \sim
\left[ \sqrt{\pi \ \kappa} \cdot \frac{\Gamma(\kappa + r - 1/2)}{\Gamma(\kappa + r)} \right] \cdot |f_0|_p, \eqno(3.4)
$$
which coincides with upper estimate (2.5b), of course  when $ p \to 2-0. $\par

\vspace{3mm}

We now give a more accurate lower bound holds for all the values $  p, \ p \in (1,2). $
Closer look closely to the function $ g_0(s). $ Consider only the case $  s \ge 1,  $ suppose in addition

$$
\kappa > 0, \hspace{6mm}  \kappa + r > 1  \eqno(3.5)
$$
and introduce a new variable

$$
Y = Y(\kappa,r) = \frac{\kappa}{\kappa + r - 1} \left[ 1 - \left( \frac{\kappa}{\kappa + 1}  \right)^{\kappa + r - 1}  \right]: \eqno(3.6)
$$

$$
\int_0^s (1 + z/\kappa)^{-\kappa - r} \frac{dz}{\sqrt{z}}\ge \int_0^1 \frac{dz}{\sqrt{z} \ (1 + z/\kappa)^{\kappa+r} } \ge
$$

$$
\int_0^1 \frac{dz}{(1 + z/\kappa)^{\kappa+r} } = Y(\kappa,r);
$$
hence

$$
|g_0|_q^q \ge Y^q \cdot \int_1^{\infty} |g_0(s)|^q \ ds \ge Y^q (\kappa,r) \cdot \frac{2}{q-2};
$$

$$
|g_0|_q:|f_0|_p \ge (p-1)^{1-1/p} \cdot (1-p/2)^{2/p - 1} \cdot Y(\kappa,r).
$$

 Thus, we proved in fact that under additional conditions (3.5)

$$
K_{\kappa,r}(p,q) \ge  (p-1)^{1-1/p} \cdot (1-p/2)^{2/p - 1} \cdot Y(\kappa,r) =
$$

$$
(p-1)^{1-1/p} \cdot (1-p/2)^{2/p - 1} \cdot
\frac{\kappa}{\kappa + r - 1} \left[ 1 - \left( \frac{\kappa}{\kappa + 1}  \right)^{\kappa + r - 1}  \right]. \eqno(3.7)
$$

 Obviously, this result remains true even in the extremal cases $ p = 1 $ and $  p = 2.$ \\

\bigskip

 \section{Weight estimates for FLT}

 \vspace{3mm}

 Let us  consider in this section the action of the fractional Laplace transform
operator $ L_{\kappa,r} $ on the {\it  weight} function $  f_{\mu}(t) = t^{\mu  - 1} \ f(t), $ where $ f(\cdot) \in L_p(R_+): $

$$
\Psi[f](s)= \Psi_{\mu,\kappa,r}[f](s):= \int_0^{\infty} t^{\mu - 1} \ (1 + ts/\kappa)^{- \kappa - r} \ f(t) \ dt = L_{\kappa,r}[f_{\mu}](s). \eqno(4.1)
$$

\vspace{3mm}

{\it  Notations and restrictions:} $ 1/\sigma := 1 + \mu - \frac{2}{p}, \ p \le Q, $

$$
 \frac{1}{Q} = \mu - \frac{1}{p}, \ p' = p/(p-1), \eqno(4.2)
$$

$$
0 < \mu < 1, \  \sigma < p' \ \Leftrightarrow 1/\mu < p \le 2/\mu,  \ (\kappa + r) \cdot \sigma + \frac{\sigma}{p'} > 1, \eqno(4.3)
$$

$$
\Theta(p) = \Theta_{\mu, \kappa,r}(p) := \left[ \kappa^{1 - \sigma/p} \ B(1 - \sigma/p', \ (\kappa + r) \cdot \sigma + \sigma/p' - 1) \right]^{1/\sigma}, \eqno(4.4)
$$
where  $ B(\alpha,\beta) $ denotes the ordinary Beta function. \par

 Notice that the relation (4.2) is necessary for the inequality of a form

$$
 \Psi_{\mu,\kappa,r}[f] \ |_Q \le  K^{\mu}_{\kappa, r}(p,Q) \cdot |f|_p,
$$
where we accept as before

$$
K^{(\mu)}_{\kappa, r}(p,Q) = \sup_{0 < |f|_p <\infty} \left[ \frac{| \ \Psi_{\mu,\kappa,r}[f] \ |_Q }{|f|_p}  \right].
$$
 This proposition may be proved by means of scaling method.\\

\vspace{3mm}

{\bf Theorem 4.1.}

$$
| \ \Psi_{\mu,\kappa,r}[f] \ |_Q \le \Theta(p) \cdot |f|_p, \eqno(4.5)
$$
{\it or equally}

$$
K^{(\mu)}_{\kappa, r}(p,Q) \le \Theta_{\mu, \kappa,r}(p).  \eqno(4.5a)
$$

\vspace{3mm}

{\bf Proof} follows immediately after simple calculations from  \cite{Okikiolu1}, p. 220-222.  \par
 Consider the following  linear operator

$$
U_{\sigma, \mu}[f](x) := \int_R |t|^{\mu - 1} \ h(x \cdot t) \ f(t) \ dt,  \eqno(4.6)
$$
where

$$
1 \le p \le Q; \ \frac{1}{Q} = \mu -  \frac{1}{p}, \  \frac{1}{\sigma } = 1 + \mu - \frac{2}{p}, \ t,x \in R. \eqno(4.6a)
$$

Denote

$$
M(p) = \left[ \int_R |t|^{- \sigma(p-1)/p } \ |h(t)|^{\sigma} \ dt   \right]^{1/\sigma}. \eqno(4.7)
$$

 It is proved in \cite{Okikiolu1}, p. 220-222 that

$$
|U_{\sigma, \mu}[f]|_Q \le M(p) \ |f|_p. \eqno(4.8)
$$

 It remains to calculate the integral $  M(p) $  in (4.7), in which we substitute

$$
  h(x) = (1 + |x|/\kappa)^{-\kappa - r}, \ f(x) = f(|x|).
$$

\vspace{3mm}

We now turn to the conclusion of the lower bound for the value  $ K^{(\mu)}_{\kappa, r}(p,Q). $ Note that the unique
critical value of the parameter $ p $ is the value $  p_0 = 2/\mu -0 $ and correspondingly $ Q = 2/\mu + 0. $
We therefore consider the following trial function

$$
f_0(x) = x^{-\mu/2} \cdot I_{(0,1)}(x);
$$
then

$$
| f_0 |_p  = \int_0^1 x^{ -p \mu/2 } dx = \frac{2}{2 - p \mu}, \ |f_0|_p = \left[ \frac{2}{2 - p \mu}   \right]^{1/p}.
$$
 Further, let us denote  for the values $  s \ge 1 $ only

$$
g_0(s) = \Psi[f_0](s) = \int_0^1 t^{\mu/2 - 1} \ (1 + ts/\kappa)^{-\kappa - r} \ dt =
$$

$$
s^{-\mu/2} \int_0^s z^{\mu/2 - 1} \ (1 + z/\kappa)^{ - \kappa - r } \ dz;
$$
then

$$
g_0(s) \ge s^{-\mu/2} \int_0^1 (1 + z/\kappa)^{ - \kappa - r } \ dz =  s^{-\mu/2} \ Y(\kappa,r);
$$

$$
|g_0|_Q \ge Y(\kappa,r) \cdot \left[ \frac{2}{\mu Q - 2} \right]^{1/Q};
$$

$$
|g_0|_Q:|f_0|_p \ge Y(\kappa,r) \cdot (1 - \mu p/2)^{2/p - \mu} \cdot (\mu p - 1)^{\mu - 1/p}.
$$

 Thus,

$$
K^{(\mu)}_{\kappa, r}(p,Q) \ge Y(\kappa,r) \cdot (1 - \mu p/2)^{2/p - \mu} \cdot (\mu p - 1)^{\mu - 1/p}, \ 1/\mu \le p \le 2/\mu.\eqno(4.9)
$$

\bigskip

\section{Generalization on the Grand Lebesgue Spaces (GLS). }

\vspace{3mm}

We recall first of all here  for reader conventions some definitions and facts from
the theory of GLS spaces.\par

Recently, see
\cite{Fiorenza1}, \cite{Fiorenza2},\cite{Ivaniec1}, \cite{Ivaniec2}, \cite{Jawerth1},
\cite{Karadzov1}, \cite{Kozatchenko1}, \cite{Liflyand1}, \cite{Ostrovsky1}, \cite{Ostrovsky2} etc.
 appear the so-called Grand Lebesgue Spaces (GLS)
 $$
 G(\psi) = G = G(\psi ; A;B);  \ A;B = \const; \ A \ge 1, \ B \le \infty
 $$
spaces consisting on all the measurable functions $ f : X \to R  $ with finite norms

$$
||f||G(\psi) \stackrel{def}{=} \sup_{p \in (A;B)} \left[\frac{|f|_p}{\psi(p)} \right]. \eqno(5.1)
$$

 Here $ \psi = \psi(p), \ p \in (A,B) $ is some continuous positive on the {\it open} interval $ (A;B) $ function such
that

$$
\inf_{p \in(A;B)} \psi(p) > 0. \eqno(5.2)
$$

We will denote
$$
\supp(\psi) \stackrel{def}{=} (A;B).
$$

The set of all such a functions with support $ \supp(\psi) = (A;B) $ will be denoted by  $  \Psi(A;B). $  \par

This spaces are rearrangement invariant; and are used, for example, in
the theory of Probability, theory of Partial Differential Equations,
 Functional Analysis, theory of Fourier series,
 Martingales, Mathematical Statistics, theory of Approximation  etc. \par

 Notice that the classical Lebesgue - Riesz spaces $ L_p $  are extremal case of Grand Lebesgue Spaces, see
 \cite{Ostrovsky2},  \cite{Ostrovsky100}. \par

 Let a function $  \xi:  R_+ \to R  $ be such that

 $$
 \exists (A,B): \ 1 \le A < B \le \infty \ \Rightarrow  \forall p \in (A,B) \ |f|_p < \infty.
 $$
Then the function $  \psi = \psi_{\xi}(p) $ may be naturally defined by the following way:

$$
\psi_{\xi}(p) := |\xi|_p, \ p \in (A,B). \eqno(5.3)
$$

\vspace{4mm}

 Let now the (measurable) function $ f: R_+ \to R, \ f \in G\psi \  $  for some $  \psi(\cdot) $ with support
 $ \supp \psi = (A,B)  $ for which

 $$
   (a,b) := (A,B) \cap (1, 2) \ne \emptyset.  \eqno(5.4)
 $$

 We define the function

 $$
\lambda(q) = \frac{q}{q-1}, \ q \in (b',a')
 $$
 and introduce a new  $ \psi \ - $ function, say $ \nu = \nu_{z} =  \nu_{z}(q)  $  as follows.

 $$
 \nu_{z}(q) =  z_{\kappa, r}(\lambda(q)) \cdot \psi(\lambda(q)), \ q \in (b',a'). \eqno(5.5)
 $$

\vspace{4mm}

{\bf Theorem 5.1.} {\it We assert  under condition (5.4):}

$$
||L_{\kappa, r} [f] ||G \nu \le 1 \cdot ||f||G\psi_{a,b},  \eqno(5.6)
$$
{\it  where the constant "1" is the best possible.} \par

\vspace{4mm}

{\bf Proof. Upper bound.} \par
 Let further  in this section $ p \in (a,b). $ We can and will suppose without loss of generality  $  ||f||G\psi_{a,b} = 1. $
Then

$$
|f|_p \le \psi_{a,b}(p), \ p \in (a,b). \eqno(5.7)
$$

 We conclude after  substituting into the  proposition of theorem 2.1

 $$
 |L_{\kappa,r}[f] |_{q} \le K_{\kappa, r}(\lambda(q) \cdot  \psi_{a,b}(\lambda(q)) = \nu_z(q), \  q \in (b',a'). \eqno(5.8)
 $$
The inequality (5.6) follows from  (5.8) after substitution $ p = \lambda(q). $ \par

\vspace{4mm}

{\bf Proof. Exactness.} \par
 The exactness of the constant "1"  in the proposition (5.8) follows immediately from  the  theorem 2.1 in the article
\cite{Ostrovsky100}. \par

\hfill $\Box$

\bigskip

\section{Concluding remarks}

\vspace{3mm}

 {\bf A. } The majority of obtained results  may be obtained for the more general operators of the form

$$
T_K[f](x) = \int_0^{\infty} K(x \cdot y) \ f(y) \ dy,  \ x \ge 0, \eqno(6.1)
$$
or more generally

$$
T_K^{(\mu-1)}[f](x) = \int_0^{\infty} x^{\mu - 1} \ K(x \cdot y) \ f(y) \ dy. \eqno(6.1a)
$$

 For instance, assume $ K(z) \ge 0 $ and define the Mellin's transform of the kernel $  K(\cdot):  $

$$
\zeta(s) = \int_0^{\infty} x^{s-1} \ K(x) \ dx.
$$
 Then

$$
T_K[f]_p^p \le \zeta(1/p) \int_0^{\infty} (x |f(x)|)^p \ \frac{dx}{x^2},
$$

$$
\int_0^{\infty} x^{p-2} \ |T_K[f](x)|^p \ dx \le \zeta^p(1 - 1/p) \ |f|_p^p,
$$
see \cite{Hardy1}, p. 256; \cite{Mitrinovic1},  p.146. \par
 Another approach. Let $ p,q = \const \ge 1, \ 1/p + 1/q = 1. $ We derive using H\"older's inequality:

 $$
 | T_K[f](x)| \le  \left[\int_0^{\infty} |K(x y)|^q \ dy  \right]^{1/q} \cdot   |f|_p =
 $$

 $$
  x^{-1/q} \cdot |K|_q \cdot |f|_p,
 $$
or equally

$$
\sup_{x > 0} \left[  x^{1/q} \ |T_K[f](x)|  \right] \le |K|_q \cdot |f|_p, \eqno(6.2)
$$
if of course $ K(\cdot) \in L_q. $ We do not suppose in (6.2) the non-negativity of the function $  K(\cdot). $ \par

\vspace{3mm}

{\bf B.} Analogously may be investigated  the so-called "multivariate" case:

$$
T_K[f](\vec{x}) = \int_{R^d} K(\vec{x} \odot \vec{y}) \ f( \vec{y}) \ d \vec{y},  \ x,y \in R^d, \eqno(6.3)
$$
where

$$
\vec{x} \odot \vec{y} = \{ x_1 y_1, x_2 y_2, \ldots, x_d y_d \}.
$$
 Needed here  $ L_q(R^d) - L_p(R^d) $  estimations for operator $ T_K[\cdot] $  may be found in the book of G.O.Okikiolu \cite{Okikiolu1},
p. 222-224, \ 341-345. \par

\end{document}